\documentclass{amsart}
\usepackage{amsfonts}
\usepackage{graphicx}
\usepackage{amscd}

\newtheorem{theorem}{Theorem}
\theoremstyle{plain}

\newtheorem{corollary}{Corollary}

\newtheorem{remark}{Remark}

\numberwithin{equation}{section}

\input{tcilatex}

\begin{document}
\title[Ostrowski Type Inequalites]{An Inequality of Ostrowski Type via
Pompeiu's Mean Value Theorem}
\author{S.S. Dragomir}
\address{School of Computer Science and Mathematics\\
Victoria University of Technology\\
PO Box 14428, MCMC 8001\\
VIC, Australia.}
\email{sever@matilda.vu.edu.au}
\urladdr{http://rgmia.vu.edu.au/SSDragomirWeb.html}
\date{February 27, 2003.}

\begin{abstract}
An inequality providing some bounds for the integral mean via Pompeiu's mean
value theorem and applications for quadrature rules and special means are
given.
\end{abstract}

\maketitle

\section{Introduction}

The following result is known in the literature as Ostrowski's inequality 
\cite{AO}.

\begin{theorem}
\label{t1}Let $f:\left[ a,b\right] \rightarrow \mathbb{R}$ be a
differentiable mapping on $\left( a,b\right) $ with the property that $%
\left\vert f^{\prime }\left( t\right) \right\vert \leq M$ for all $t\in
\left( a,b\right) .$ Then 
\begin{equation}
\left\vert f\left( x\right) -\frac{1}{b-a}\int_{a}^{b}f\left( t\right)
dt\right\vert \leq \left[ \frac{1}{4}+\left( \frac{x-\frac{a+b}{2}}{b-a}%
\right) ^{2}\right] \left( b-a\right) M,  \label{1.1}
\end{equation}%
for all $x\in \left[ a,b\right] .$ The constant $\frac{1}{4}$ is best
possible in the sense that it cannot be replaced by a smaller constant.
\end{theorem}

In \cite{SSD1}, the author has proved the following Ostrowski type
inequality.

\begin{theorem}
\label{t2}Let $f:\left[ a,b\right] \rightarrow \mathbb{R}$ be continuous on $%
\left[ a,b\right] $ with $a>0$ and differentiable on $\left( a,b\right) .$
Let $p\in \mathbb{R}\backslash \left\{ 0\right\} $ and assume that 
\begin{equation*}
K_{p}\left( f^{\prime }\right) :=\sup\limits_{u\in \left( a,b\right)
}\left\{ u^{1-p}\left| f^{\prime }\left( u\right) \right| \right\} <\infty .
\end{equation*}
Then we have the inequality 
\begin{multline}
\left| f\left( x\right) -\frac{1}{b-a}\int_{a}^{b}f\left( t\right) dt\right|
\leq \frac{K_{p}\left( f^{\prime }\right) }{\left| p\right| \left(
b-a\right) }  \label{1.2} \\
\times \left\{ 
\begin{array}{ll}
2x^{p}\left( x-A\right) +\left( b-x\right) L_{p}^{p}\left( b,x\right)
-\left( x-a\right) L_{p}^{p}\left( x,a\right) & \text{if }p\in \left(
0,\infty \right) ; \\ 
&  \\ 
\left( x-a\right) L_{p}^{p}\left( x,a\right) -\left( b-x\right)
L_{p}^{p}\left( b,x\right) -2x^{p}\left( x-A\right) & \text{if }p\in \left(
-\infty ,-1\right) \cup \left( -1,0\right) \\ 
&  \\ 
\left( x-a\right) L^{-1}\left( x,a\right) -\left( b-x\right) L^{-1}\left(
b,x\right) -\frac{2}{x}\left( x-A\right) & \text{if }p=-1,%
\end{array}
\right.
\end{multline}
for any $x\in \left( a,b\right) ,$ where for $a\neq b,$%
\begin{equation*}
A=A\left( a,b\right) :=\frac{a+b}{2},\text{ \hspace{0.05in}is the arithmetic
mean,}
\end{equation*}
\begin{equation*}
L_{p}=L_{p}\left( a,b\right) =\left[ \frac{b^{p+1}-a^{p+1}}{\left(
p+1\right) \left( b-a\right) }\right] ^{\frac{1}{p}},\;\text{is the }p-\text{%
logarithmic mean}\;p\in \mathbb{R}\backslash \left\{ -1,0\right\} ,
\end{equation*}
and 
\begin{equation*}
L=L\left( a,b\right) :=\frac{b-a}{\ln b-\ln a}\text{ \hspace{0.05in} is the
logarithmic mean.}
\end{equation*}
\end{theorem}

Another result of this type obtained in the same paper is:

\begin{theorem}
\label{t3}Let $f:\left[ a,b\right] \rightarrow \mathbb{R}$ be continuous on $%
\left[ a,b\right] $ (with $a>0$) and differentiable on $\left( a,b\right) .$
If 
\begin{equation*}
P\left( f^{\prime }\right) :=\sup\limits_{u\in \left( a,b\right) }\left|
uf^{\prime }\left( x\right) \right| <\infty ,
\end{equation*}
then we have the inequality 
\begin{equation}
\left| f\left( x\right) -\frac{1}{b-a}\int_{a}^{b}f\left( t\right) dt\right|
\leq \frac{P\left( f^{\prime }\right) }{b-a}\left[ \ln \left[ \frac{\left[
I\left( x,b\right) \right] ^{b-x}}{\left[ I\left( a,x\right) \right] ^{x-a}}%
\right] +2\left( x-A\right) \ln x\right]  \label{1.3}
\end{equation}
for any $x\in \left( a,b\right) ,$ where for $a\neq b$%
\begin{equation*}
I=I\left( a,b\right) :=\frac{1}{e}\left( \frac{b^{b}}{a^{a}}\right) ^{\frac{1%
}{b-a}},\text{ \hspace{0.05in} is the identric mean.}
\end{equation*}
\end{theorem}

If some local information around the point $x\in \left( a,b\right) $ is
available, then we may state the following result as well \cite{SSD1}.

\begin{theorem}
\label{t4}Let $f:\left[ a,b\right] \rightarrow \mathbb{R}$ be continuous on $%
\left[ a,b\right] $ and differentiable on $\left( a,b\right) .$ Let $p\in
\left( 0,\infty \right) $ and assume, for a given $x\in \left( a,b\right) ,$
we have that 
\begin{equation*}
M_{p}\left( x\right) :=\sup\limits_{u\in \left( a,b\right) }\left\{
\left\vert x-u\right\vert ^{1-p}\left\vert f^{\prime }\left( u\right)
\right\vert \right\} <\infty .
\end{equation*}%
Then we have the inequality 
\begin{multline}
\left\vert f\left( x\right) -\frac{1}{b-a}\int_{a}^{b}f\left( t\right)
dt\right\vert  \label{1.4} \\
\leq \frac{1}{p\left( p+1\right) \left( b-a\right) }\left[ \left( x-a\right)
^{p+1}+\left( b-x\right) ^{p+1}\right] M_{p}\left( x\right) .
\end{multline}
\end{theorem}

For recent results in connection to Ostrowski's inequality see the papers 
\cite{GA1},\cite{GA2} and the monograph \cite{SSD2}.

The main aim of this paper is to provide some complementary results, where
instead of using Cauchy mean value theorem, we use Pompeiu mean Value
Theorem to evaluate the integral mean of an absolutely continuous function.
Applications for quadrature rules and particular instances of functions are
given as well.

\section{Pompeiu's Mean Value Theorem}

In 1946, Pompeiu \cite{POM} derived a variant of Lagrange's mean value
theorem, now known as \textit{Pompeiu's mean value theorem} (see also \cite[%
p. 83]{SR}).

\begin{theorem}
\label{t2.1}For every real valued function $f$ differentiable on an interval 
$\left[ a,b\right] $ not containing $0$ and for all pairs $x_{1}\neq x_{2}$
in $\left[ a,b\right] ,$ there exists a point $\xi $ in $\left(
x_{1},x_{2}\right) $ such that 
\begin{equation}
\frac{x_{1}f\left( x_{2}\right) -x_{2}f\left( x_{1}\right) }{x_{1}-x_{2}}%
=f\left( \xi \right) -\xi f^{\prime }\left( \xi \right) .  \label{2.1}
\end{equation}
\end{theorem}

\begin{proof}
Define a real valued function $F$ on the interval $\left[ \frac{1}{b},\frac{1%
}{a}\right] $ by 
\begin{equation}
F\left( t\right) =tf\left( \frac{1}{t}\right) .  \label{2.2}
\end{equation}%
Since $f$ is differentiable on $\left( \frac{1}{b},\frac{1}{a}\right) $ and 
\begin{equation}
F^{\prime }\left( t\right) =f\left( \frac{1}{t}\right) -\frac{1}{t}f^{\prime
}\left( \frac{1}{t}\right) ,  \label{2.3}
\end{equation}%
then applying the mean value theorem to $F$ on the interval $\left[ x,y%
\right] \subset \left[ \frac{1}{b},\frac{1}{a}\right] $ we get 
\begin{equation}
\frac{F\left( x\right) -F\left( y\right) }{x-y}=F^{\prime }\left( \eta
\right)  \label{2.4}
\end{equation}%
for some $\eta \in \left( x,y\right) .$

Let $x_{2}=\frac{1}{x}$, $x_{1}=\frac{1}{y}$ and $\xi =\frac{1}{\eta }.$
Then, since $\eta \in \left( x,y\right) ,$ we have 
\begin{equation*}
x_{1}<\xi <x_{2}.
\end{equation*}
Now, using (\ref{2.2}) and (\ref{2.3}) on (\ref{2.4}), we have 
\begin{equation*}
\frac{xf\left( \frac{1}{x}\right) -yf\left( \frac{1}{y}\right) }{x-y}%
=f\left( \frac{1}{\eta }\right) -\frac{1}{\eta }f^{\prime }\left( \frac{1}{%
\eta }\right) ,
\end{equation*}
that is 
\begin{equation*}
\frac{x_{1}f\left( x_{2}\right) -x_{2}f\left( x_{1}\right) }{x_{1}-x_{2}}%
=f\left( \xi \right) -\xi f^{\prime }\left( \xi \right) .
\end{equation*}
This completes the proof of the theorem.
\end{proof}

\begin{remark}
Following \cite[p. 84 -- 85]{SR}, we will mention here a geometrical
interpretation of Pompeiu's theorem. The equation of the secant line joining
the points $\left( x_{1},f\left( x_{1}\right) \right) $ and $\left(
x_{2},f\left( x_{2}\right) \right) $ is given by 
\begin{equation*}
y=f\left( x_{1}\right) +\frac{f\left( x_{2}\right) -f\left( x_{1}\right) }{%
x_{2}-x_{1}}\left( x-x_{1}\right) .
\end{equation*}%
This line intersects the $y-$axis at the point $\left( 0,y\right) ,$ where $%
y $ is 
\begin{align*}
y& =f\left( x_{1}\right) +\frac{f\left( x_{2}\right) -f\left( x_{1}\right) }{%
x_{2}-x_{1}}\left( 0-x_{1}\right) \\
& =\frac{x_{1}f\left( x_{2}\right) -x_{2}f\left( x_{1}\right) }{x_{1}-x_{2}}.
\end{align*}%
The equation of the tangent line at the point $\left( \xi ,f\left( \xi
\right) \right) $ is 
\begin{equation*}
y=\left( x-\xi \right) f^{\prime }\left( \xi \right) +f\left( \xi \right) .
\end{equation*}%
The tangent line intersects the $y-$axis at the point $\left( 0,y\right) ,$
where 
\begin{equation*}
y=-\xi f^{\prime }\left( \xi \right) +f\left( \xi \right) .
\end{equation*}%
Hence, the geometric meaning of Pompeiu's mean value theorem is that the
tangent of the point $\left( \xi ,f\left( \xi \right) \right) $ intersects
on the $y-$axis at the same point as the secant line connecting the points $%
\left( x_{1},f\left( x_{1}\right) \right) $ and $\left( x_{2},f\left(
x_{2}\right) \right) .$
\end{remark}

\section{Evaluating the Integral Mean}

The following result holds.

\begin{theorem}
\label{t3.1}Let $f:\left[ a,b\right] \rightarrow \mathbb{R}$ be continuous
on $\left[ a,b\right] $ and differentiable on $\left( a,b\right) $ with $%
\left[ a,b\right] $ not containing $0.$ Then for any $x\in \left[ a,b\right]
,$ we have the inequality 
\begin{multline}
\qquad \left| \frac{a+b}{2}\cdot \frac{f\left( x\right) }{x}-\frac{1}{b-a}%
\int_{a}^{b}f\left( t\right) dt\right|  \label{3.1} \\
\leq \frac{b-a}{\left| x\right| }\left[ \frac{1}{4}+\left( \frac{x-\frac{a+b%
}{2}}{b-a}\right) ^{2}\right] \left\| f-\ell f^{\prime }\right\| _{\infty
},\qquad
\end{multline}
where $\ell \left( t\right) =t,$ $t\in \left[ a,b\right] .$

The constant $\frac{1}{4}$ is sharp in the sense that it cannot be replaced
by a smaller constant.
\end{theorem}

\begin{proof}
Applying Pompeiu's mean value theorem, for any $x,t\in \left[ a,b\right] ,$
there is a $\xi $ between $x$ and $t$ such that 
\begin{equation*}
tf\left( x\right) -xf\left( t\right) =\left[ f\left( \xi \right) -\xi
f^{\prime }\left( \xi \right) \right] \left( t-x\right)
\end{equation*}%
giving 
\begin{align}
\left\vert tf\left( x\right) -xf\left( t\right) \right\vert & \leq
\sup\limits_{\xi \in \left[ a,b\right] }\left\vert f\left( \xi \right) -\xi
f^{\prime }\left( \xi \right) \right\vert \left\vert x-t\right\vert
\label{3.2} \\
& =\left\Vert f-\ell f^{\prime }\right\Vert _{\infty }\left\vert
x-t\right\vert  \notag
\end{align}%
for any $t,x\in \left[ a,b\right] .$

Integrating over $t\in \left[ a,b\right] ,$ we get 
\begin{align}
\left\vert f\left( x\right) \int_{a}^{b}tdt-x\int_{a}^{b}f\left( t\right)
dt\right\vert & \leq \left\Vert f-\ell f^{\prime }\right\Vert _{\infty
}\int_{a}^{b}\left\vert x-t\right\vert dt  \label{3.3} \\
& =\left\Vert f-\ell f^{\prime }\right\Vert _{\infty }\left[ \frac{\left(
x-a\right) ^{2}+\left( b-x\right) ^{2}}{2}\right]  \notag \\
& =\left\Vert f-\ell f^{\prime }\right\Vert _{\infty }\left[ \frac{1}{4}%
\left( b-a\right) ^{2}+\left( x-\frac{a+b}{2}\right) ^{2}\right]  \notag
\end{align}%
and since $\int_{a}^{b}tdt=\frac{b^{2}-a^{2}}{2},$ we deduce from (\ref{3.3}%
) the desired result (\ref{3.1}).

Now, assume that (\ref{3.2}) holds with a constant $k>0,$ i.e., 
\begin{multline}
\qquad \left| \frac{a+b}{2}\cdot \frac{f\left( x\right) }{x}-\frac{1}{b-a}%
\int_{a}^{b}f\left( t\right) dt\right|  \label{3.4} \\
\leq \frac{b-a}{\left| x\right| }\left[ k+\left( \frac{x-\frac{a+b}{2}}{b-a}%
\right) ^{2}\right] \left\| f-\ell f^{\prime }\right\| _{\infty },\qquad
\end{multline}
for any $x\in \left[ a,b\right] .$

Consider $f:\left[ a,b\right] \rightarrow \mathbb{R}$, $f\left( t\right)
=\alpha t+\beta ;$ $\alpha ,\beta \neq 0.$ Then 
\begin{eqnarray*}
\left\Vert f-\ell f^{\prime }\right\Vert _{\infty } &=&\left\vert \beta
\right\vert , \\
\frac{1}{b-a}\int_{a}^{b}f\left( t\right) dt &=&\frac{a+b}{2}\cdot \alpha
+\beta ,
\end{eqnarray*}%
and by (\ref{3.4}) we deduce 
\begin{equation*}
\left\vert \frac{a+b}{2}\left( \alpha +\frac{\beta }{x}\right) -\left( \frac{%
a+b}{2}\alpha +\beta \right) \right\vert \leq \frac{b-a}{\left\vert
x\right\vert }\left[ k+\left( \frac{x-\frac{a+b}{2}}{b-a}\right) ^{2}\right]
\left\vert \beta \right\vert
\end{equation*}%
giving 
\begin{equation}
\left\vert \frac{a+b}{2}-x\right\vert \leq \left( b-a\right) k+\left( \frac{%
x-\frac{a+b}{2}}{b-a}\right) ^{2}  \label{3.5}
\end{equation}%
for any $x\in \left[ a,b\right] .$

If in (\ref{3.5}) we let $x=a$ or $x=b$, we deduce $k\geq \frac{1}{4}$, and
the sharpness of the constant is proved.
\end{proof}

The following interesting particular case holds.

\begin{corollary}
\label{c3.2}With the assumptions in Theorem \ref{t3.1}, we have 
\begin{equation}
\left| f\left( \frac{a+b}{2}\right) -\frac{1}{b-a}\int_{a}^{b}f\left(
t\right) dt\right| \leq \frac{\left( b-a\right) }{2\left| a+b\right| }%
\left\| f-\ell f^{\prime }\right\| _{\infty }.  \label{3.6}
\end{equation}
\end{corollary}

\section{The Case of Weighted Integrals}

We will consider now the weighted integral case.

\begin{theorem}
\label{t4.1}Let $f:\left[ a,b\right] \rightarrow \mathbb{R}$ be continuous
on $\left[ a,b\right] $ and differentiable on $\left( a,b\right) $ with $%
\left[ a,b\right] $ not containing $0.$ If $w:\left[ a,b\right] \rightarrow 
\mathbb{R}$ is nonnegative integrable on $\left[ a,b\right] ,$ then for each 
$x\in \left[ a,b\right] ,$ we have the inequality: 
\begin{multline}
\left| \int_{a}^{b}f\left( t\right) w\left( t\right) dt-\frac{f\left(
x\right) }{x}\int_{a}^{b}tw\left( t\right) dt\right|  \label{4.1} \\
\leq \left\| f-\ell f^{\prime }\right\| _{\infty }\left[ \func{sgn}\left(
x\right) \left( \int_{a}^{x}w\left( t\right) dt-\int_{x}^{b}w\left( t\right)
dt\right) \right. \\
+\left. \frac{1}{\left| x\right| }\left( \int_{x}^{b}tw\left( t\right)
dt-\int_{a}^{x}tw\left( t\right) dt\right) \right]
\end{multline}
\end{theorem}

\begin{proof}
Using the inequality (\ref{3.2}), we have 
\begin{align}
& \left| f\left( x\right) \int_{a}^{b}tw\left( t\right)
dt-x\int_{a}^{b}f\left( t\right) w\left( t\right) dt\right|  \label{4.2} \\
& \leq \left\| f-\ell f^{\prime }\right\| _{\infty }\int_{a}^{b}w\left(
t\right) \left| x-t\right| dt  \notag \\
& =\left\| f-\ell f^{\prime }\right\| _{\infty }\left[ \int_{a}^{x}w\left(
t\right) \left( x-t\right) dt+\int_{x}^{b}w\left( t\right) \left( t-x\right)
dt\right]  \notag \\
& =\left\| f-\ell f^{\prime }\right\| _{\infty }\left[ x\int_{a}^{x}w\left(
t\right) dt-\int_{a}^{x}tw\left( t\right) dt+\int_{x}^{b}tw\left( t\right)
dt-x\int_{x}^{b}w\left( t\right) dt\right]  \notag \\
& =\left\| f-\ell f^{\prime }\right\| _{\infty }\left[ x\left(
\int_{a}^{x}w\left( t\right) dt-\int_{x}^{b}w\left( t\right) dt\right)
+\int_{x}^{b}tw\left( t\right) dt-\int_{a}^{x}tw\left( t\right) dt\right] 
\notag
\end{align}
from where we get the desired inequality (\ref{4.1}).
\end{proof}

Now, if we assume that $0<a<b,$ then 
\begin{equation}
a\leq \frac{\int_{a}^{b}tw\left( t\right) dt}{\int_{a}^{b}w\left( t\right) dt%
}\leq b  \label{4.3}
\end{equation}
provided $\int_{a}^{b}w\left( t\right) dt>0.$

With this extra hypothesis, we may state the following corollary.

\begin{corollary}
\label{c4.2}With the above assumptions, we have 
\begin{multline}
\left| f\left( \frac{\int_{a}^{b}tw\left( t\right) dt}{\int_{a}^{b}w\left(
t\right) dt}\right) -\frac{1}{\int_{a}^{b}w\left( t\right) dt}%
\int_{a}^{b}f\left( t\right) w\left( t\right) dt\right|  \label{4.4} \\
\leq \left\| f-\ell f^{\prime }\right\| _{\infty }\left[ \frac{%
\int_{a}^{x}w\left( t\right) dt-\int_{x}^{b}w\left( t\right) dt}{%
\int_{a}^{b}w\left( t\right) dt}+\frac{\int_{x}^{b}w\left( t\right)
tdt-\int_{a}^{x}tw\left( t\right) dt}{\int_{a}^{b}tw\left( t\right) dt}%
\right]
\end{multline}
\end{corollary}

\section{A Quadrature Formula}

We assume in the following that $0<a<b.$

Consider the division of the interval $\left[ a,b\right] $ given by%
\begin{equation*}
I_{n}:a=x_{0}<x_{1}<\cdots <x_{n-1}<x_{n}=b,
\end{equation*}%
and $\xi _{i}\in \left[ x_{i},x_{i+1}\right] $, $i=0,\dots ,n-1$ a sequence
of intermediate points. Define the quadrature 
\begin{align}
S_{n}\left( f,I_{n},\mathbf{\xi }\right) & :=\sum_{i=0}^{n-1}\frac{f\left(
\xi _{i}\right) }{\xi _{i}}\cdot \frac{x_{i+1}^{2}-x_{i}^{2}}{2}  \label{5.1}
\\
& =\sum_{i=0}^{n-1}\frac{f\left( \xi _{i}\right) }{\xi _{i}}\cdot \frac{%
x_{i}+x_{i+1}}{2}\cdot h_{i},  \notag
\end{align}%
where $h_{i}:=x_{i+1}-x_{i},$ $i=0,\dots ,n-1.$

The following result concerning the estimate of the remainder in
approximating the integral $\int_{a}^{b}f\left( t\right) dt$ by the use of $%
S_{n}\left( f,I_{n},\mathbf{\xi }\right) $ holds.

\begin{theorem}
\label{t5.1}Assume that $f:\left[ a,b\right] \rightarrow \mathbb{R}$ is
continuous on $\left[ a,b\right] $ and differentiable on $\left( a,b\right)
. $ Then we have the representation 
\begin{equation}
\int_{a}^{b}f\left( t\right) dt=S_{n}\left( f,I_{n},\mathbf{\xi }\right)
+R_{n}\left( f,I_{n},\mathbf{\xi }\right)  \label{5.2}
\end{equation}%
where $S_{n}\left( f,I_{n},\mathbf{\xi }\right) $ is as defined in (\ref{5.1}%
), and the remainder $R_{n}\left( f,I_{n},\mathbf{\xi }\right) $ satisfies
the estimate 
\begin{align}
\left\vert R_{n}\left( f,I_{n},\mathbf{\xi }\right) \right\vert & \leq
\left\Vert f-\ell f^{\prime }\right\Vert _{\infty }\sum_{i=0}^{n-1}\frac{%
h_{i}^{2}}{\xi _{i}}\left[ \frac{1}{4}+\left\vert \frac{\xi _{i}-\frac{%
x_{i}+x_{i+1}}{2}}{h_{i}}\right\vert ^{2}\right]  \label{5.3} \\
& \leq \frac{1}{2}\left\Vert f-\ell f^{\prime }\right\Vert _{\infty
}\sum_{i=0}^{n-1}\frac{h_{i}^{2}}{\xi _{i}}\leq \frac{\left\Vert f-\ell
f^{\prime }\right\Vert _{\infty }}{2a}\sum_{i=0}^{n-1}h_{i}^{2}.  \notag
\end{align}
\end{theorem}

\begin{proof}
Apply Theorem \ref{t3.1} on the interval $\left[ x_{i},x_{i+1}\right] $ for
the intermediate points $\xi _{i}$ to obtain 
\begin{align}
& \left\vert \int_{x_{i}}^{x_{i+1}}f\left( t\right) dt-\frac{f\left( \xi
_{i}\right) }{\xi _{i}}\cdot \frac{x_{i}+x_{i+1}}{2}\cdot h_{i}\right\vert
\label{5.4} \\
& \leq \frac{1}{\xi _{i}}h_{i}^{2}\left[ \frac{1}{4}+\left( \frac{\xi _{i}-%
\frac{x_{i}+x_{i+1}}{2}}{h_{i}}\right) ^{2}\right] \left\Vert f-\ell
f^{\prime }\right\Vert _{\infty }  \notag \\
& \leq \frac{1}{2\xi _{i}}h_{i}^{2}\left\Vert f-\ell f^{\prime }\right\Vert
_{\infty }\leq \frac{1}{2a}h_{i}^{2}\left\Vert f-\ell f^{\prime }\right\Vert
_{\infty }  \notag
\end{align}%
for each $i\in \left\{ 0,\dots ,n-1\right\} .$

Summing over $i$ from $1$ to $n-1$ and using the generalised triangle
inequality, we deduce the desired estimate (\ref{5.3}).
\end{proof}

Now, if we consider the mid-point rule (i.e., we choose $\xi _{i}=\frac{%
x_{i}+x_{i+1}}{2}$ above, $i\in \left\{ 0,\dots ,n-1\right\} $) 
\begin{equation*}
M_{n}\left( f,I_{n}\right) :=\sum_{i=0}^{n-1}f\left( \frac{x_{i}+x_{i+1}}{2}%
\right) h_{i},
\end{equation*}%
then, by Corollary \ref{c3.2}, we may state the following result as well.

\begin{corollary}
\label{c5.2}With the assumptions of Theorem \ref{t5.1}, we have 
\begin{equation}
\int_{a}^{b}f\left( t\right) dt=M_{n}\left( f,I_{n}\right) +R_{n}\left(
f,I_{n}\right) ,  \label{5.5}
\end{equation}
where the remainder satisfies the estimate: 
\begin{align}
\left| R_{n}\left( f,I_{n}\right) \right| & \leq \frac{\left\| f-\ell
f^{\prime }\right\| _{\infty }}{2}\sum_{i=0}^{n-1}\frac{h_{i}^{2}}{%
x_{i}+x_{i+1}}  \label{5.6} \\
& \leq \frac{\left\| f-\ell f^{\prime }\right\| _{\infty }}{4a}%
\sum_{i=0}^{n-1}h_{i}^{2}.  \notag
\end{align}
\end{corollary}

\section{Applications for Special Means}

For $0<a<b,$ let us consider the means 
\begin{align*}
A& =A\left( a,b\right) :=\frac{a+b}{2}, \\
G& =G\left( a,b\right) :=\sqrt{a\cdot b}, \\
H& =H\left( a,b\right) :=\frac{2}{\frac{1}{a}+\frac{1}{b}}, \\
L& =L\left( a,b\right) :=\frac{b-a}{\ln b-\ln a}\text{ \hspace{0.05in}
(logarithmic mean),} \\
I& =I\left( a,b\right) :=\frac{1}{e}\left( \frac{b^{b}}{a^{a}}\right) ^{%
\frac{1}{b-a}}\text{ \hspace{0.05in} (identric mean)}
\end{align*}%
and the $p-$logarithmic mean 
\begin{equation*}
L_{p}=L_{p}\left( a,b\right) =\left[ \frac{b^{p+1}-a^{p+1}}{\left(
p+1\right) \left( b-a\right) }\right] ^{\frac{1}{p}},\;\;p\in \mathbb{R}%
\backslash \left\{ -1,0\right\} .
\end{equation*}%
It is well known that 
\begin{equation*}
H\leq G\leq L\leq I\leq A,
\end{equation*}%
and, denoting $L_{0}:=I$ and $L_{-1}=L,$ the function $\mathbb{R}\ni
p\mapsto L_{p}\in \mathbb{R}$ is monotonic increasing.

In the following we will use the following inequality obtained in Corollary %
\ref{c3.2}, 
\begin{equation}
\left| f\left( \frac{a+b}{2}\right) -\frac{1}{b-a}\int_{a}^{b}f\left(
t\right) dt\right| \leq \frac{\left( b-a\right) }{2\left( a+b\right) }%
\left\| f-\ell f^{\prime }\right\| _{\infty },  \label{6.1}
\end{equation}
provided $0<a<b.$

\begin{enumerate}
\item Consider the function $f:\left[ a,b\right] \subset \left( 0,\infty
\right) \rightarrow \mathbb{R}$, $f\left( t\right) =t^{p},$ $p\in \mathbb{R}%
\backslash \left\{ -1,0\right\} .$ Then 
\begin{align*}
f\left( \frac{a+b}{2}\right) & =\left[ A\left( a,b\right) \right] ^{p}, \\
\frac{1}{b-a}\int_{a}^{b}f\left( t\right) dt& =L_{p}^{p}\left( a,b\right) ,
\\
\left\| f-\ell f^{\prime }\right\| _{\left[ a,b\right] ,\infty }& =\left\{ 
\begin{array}{ll}
\left( 1-p\right) a^{p} & \text{if \hspace{0.05in}}p\in \left( -\infty
,0\right) \backslash \left\{ -1\right\} , \\ 
&  \\ 
\left| 1-p\right| b^{p} & \text{if \hspace{0.05in}}p\in \left( 0,1\right)
\cup \left( 1,\infty \right) .%
\end{array}
\right.
\end{align*}
Consequently, by (\ref{6.1}) we deduce 
\begin{multline}
\left| A^{p}\left( a,b\right) -L_{p}^{p}\left( a,b\right) \right|
\label{6.2} \\
\leq \frac{1}{4}\times \left\{ 
\begin{array}{ll}
\dfrac{\left( 1-p\right) a^{p}\left( b-a\right) }{A\left( a,b\right) } & 
\text{if \hspace{0.05in}}p\in \left( -\infty ,0\right) \backslash \left\{
-1\right\} , \\ 
&  \\ 
\dfrac{\left| 1-p\right| b^{p}\left( b-a\right) }{A\left( a,b\right) } & 
\text{if \hspace{0.05in}}p\in \left( 0,1\right) \cup \left( 1,\infty \right)
.%
\end{array}
\right.
\end{multline}

\item Consider the function $f:\left[ a,b\right] \subset \left( 0,\infty
\right) \rightarrow \mathbb{R}$, $f\left( t\right) =\frac{1}{t}.$ Then 
\begin{align*}
f\left( \frac{a+b}{2}\right) & =\frac{1}{A\left( a,b\right) }, \\
\frac{1}{b-a}\int_{a}^{b}f\left( t\right) dt& =\frac{1}{L\left( a,b\right) },
\\
\left\| f-\ell f^{\prime }\right\| _{\left[ a,b\right] ,\infty }& =\frac{2}{a%
}.
\end{align*}
Consequently, by (\ref{6.1}) we deduce 
\begin{equation}
0\leq A\left( a,b\right) -L_{p}\left( a,b\right) \leq \frac{b-a}{2a}L\left(
a,b\right) .  \label{6.3}
\end{equation}

\item Consider the function $f:\left[ a,b\right] \subset \left( 0,\infty
\right) \rightarrow \mathbb{R}$, $f\left( t\right) =\ln t.$ Then 
\begin{align*}
f\left( \frac{a+b}{2}\right) & =\ln \left[ A\left( a,b\right) \right] , \\
\frac{1}{b-a}\int_{a}^{b}f\left( t\right) dt& =\ln \left[ I\left( a,b\right) %
\right] , \\
\left\Vert f-\ell f^{\prime }\right\Vert _{\left[ a,b\right] ,\infty }&
=\max \left\{ \left\vert \ln \left( \frac{a}{e}\right) \right\vert
,\left\vert \ln \left( \frac{b}{e}\right) \right\vert \right\} .
\end{align*}%
Consequently, by (\ref{6.1}) we deduce 
\begin{equation}
1\leq \frac{A\left( a,b\right) }{I\left( a,b\right) }\leq \exp \left\{ \frac{%
b-a}{4A\left( a,b\right) }\max \left\{ \left\vert \ln \left( \frac{a}{e}%
\right) \right\vert ,\left\vert \ln \left( \frac{b}{e}\right) \right\vert
\right\} \right\} .  \label{6.4}
\end{equation}
\end{enumerate}

\end{document}